\documentclass[10pt]{article}
\usepackage{amsmath,amssymb,amsbsy,graphicx,fullpage}

\let\epsilon\varepsilon
\let\phi\varphi

\def\N{\mathbb N}

\def\E{{\bf E}}
\def\C{{\mathcal C}}

\newtheorem{theorem}{Theorem}

\newtheorem{corollary}{Corollary}

\begin{document}

\title{An impossibility result for process discrimination.}
\author{
Daniil Ryabko}
\date{}
\maketitle
\begin{abstract}
Two series of binary observations $x_1,x_1,\dots$ and $y_1,y_2,\dots$  are presented: at each time $n\in\N$ we are given $x_n$ and $y_n$.
It is assumed that the sequences are generated independently of each other by two stochastic processes. 
 We are interested in the question of whether the sequences represent a typical realization of two different  processes or of the same one.
We demonstrate that this is impossible to decide in the case when the processes are  $B$-processes.
It follows that discrimination is 
impossible for the set of all (finite-valued) stationary ergodic processes in general.
This result means that  every discrimination procedure is bound
to  err with non-negligible frequency when presented with sequences from some of such processes. 
 It contrasts earlier positive results on $B$-processes, in particular
those showing that there are consistent $\bar d$-distance estimates for this class of processes. 
\end{abstract}

\section{Introduction}
Given two series of observations we wish to decide whether they were generated by the same process or by different ones.
The question is relatively simple when the time series are generated by a source of independent identically distributed
outcomes. It is far less clear how to solve the problem for more general cases, such as the case of stationary ergodic
time series. In this work we demonstrate that the question is impossible to decide even in the weakest asymptotic sense,
for a wide class of processes, which is a subset of the set of all stationary ergodic processes.

More formally, two series of binary observations $x_1,x_1,\dots$ and $y_1,y_2,\dots$  are presented  sequentially.
A {\em discrimination procedure} $D$ is a family of mappings $D_n: X^n\times X^n\rightarrow\{0,1\}$, $n\in\N$, that 
maps a pair of samples $(x_1,\dots,x_n)$, $(y_1,\dots,y_n)$ into a binary (``yes'' or ``no'') answer: the samples are generated by 
different distributions, or they are generated by the same distribution. 

A discrimination procedure $D$ is {\em asymptotically correct for a set $\C$ of process distributions} if for any two  distributions
$\rho_x,\rho_y\in\C$ independently generating the sequences $x_1,x_2,\dots$ and $y_1,y_2,\dots$ correspondingly the expected output converges to the correct answer: the following limit exists and the equality holds
$$
\lim_{n\rightarrow\infty}\E  D_n((x_1,\dots,x_n), (y_1,\dots,y_n))=\left\{\begin{array}{ll}0 & \text{ if $\rho_x=\rho_y$}\\ 1 &\text{ otherwise }\end{array}\right..
$$
Note that one can consider other notions of asymptotic correctness, for  example one can require the output to stabilize on the correct answer with 
probability~1. The notion of correctness that we consider is perhaps one of the weakest.
Clearly, asymptotically correct discriminating procedures exist for many classes of processes, for example for the class of all i.i.d. processes,
 or various parametric families, 
see e.g. \cite{cs,leh}; some realted positive results on hypothesis testing for stationary ergodic process can be found in \cite{rg, rr}.

We will show that asymptotically correct discrimination procedures do not exist for 
the class 
of $B$-processes,
or for the class of all stationary ergodic processes. 
 This result for $B$-processes  is  interesting  in view of some previously established results; thus, in \cite{ow,os} it is  shown that consistent
estimates of $\bar d$-distance  for $B$-processes (see definitions below) exist, while  it is impossible to estimate this distance outside this class (i.e. in general 
for stationary ergodic processes). 
So, our result  demonstrates that discrimination is harder than distance estimation. 
The distinction between these problems becomes very apparent in  view of the positive results of~\cite{rr}, which show that consistent change point estimates 
and process classification procedures exist for the class of  stationary ergodic processes.
The result of the present work also complements earlier 
negative results on $B$-processes, such as \cite{sh2} that shows that upper and lower divergence rates need not be the same for $B$-processes, 
and on stationary ergodic processes, such as \cite{br, gml, nob1, mor2}, that establish negative results concerning prediction, density estimation, and 
testing properties of processes. 
 It is worth noting that $B$-processes are of particular importance  for information theory, in particular, since they are what can be obtained by stationary codings of 
memoryless processes \cite{or2,sh}.

Next we briefly introduce the notation.
We are considering stationary ergodic processes  (time series), defined as probability distributions  on the set of one-way infinite 
sequences $A^\infty$, where $A=\{0,1\}$.
We will also consider stationary ergodic Markov chains on a countable set of states; for now let the set of states be $\N$.
Any function $f:\N\rightarrow A$ mapping the set of states to $A$, together with a stationary ergodic Markov chain $m$ defines a stationary ergodic
binary-valued process, whose value on  each time step is the value of $f$ applied to the current state of $m$. 

For two finite-valued stationary processes $\rho_x$ and $\rho_y$  the {\em $\bar d$-distance} $\bar d(\rho_x,\rho_y)$ is said to be less than $\epsilon$ 
 if there exists a single stationary process $\nu_{xy}$ on pairs $(x_n,y_n)$, $n\in\N$, such that $x_n$, $n\in\N$ are distributed according 
to $\rho_x$ and $y_n$ are distributed according to $\rho_y$ while 
\begin{equation}\label{eq:dbar} \nu_{xy}(x_1\ne y_1)\le \epsilon.\end{equation}
The infimum of the $\epsilon$'s for which a coupling can be found such that (\ref{eq:dbar}) is satisfied is taken to be the $\bar d$-distance
between $\rho_x$ and $\rho_y$. 
A process is called a {\em  $B$-process} (or a Bernoulli process) if it is in the $\bar d$-closure of the set of all aperiodic stationary ergodic $k$-step Markov processes, where $k\in\N$.
For more information on $\bar d$-distance and $B$-processes 
the reader is referred to \cite{ow,orns}.

\section{Main results}

The  main result of this work is the following theorem; the construction used in the proof is based on the same ideas as the construction used  in \cite{br} (see also  \cite{gml}) to demonstrate that  consistent prediction for stationary ergodic processes is impossible. 

\begin{theorem}\label{th} There is no asymptotically correct discrimination procedure for  the class of $B$-processes.
\end{theorem}
Since the class of $B$-processes is a subset of the class of all stationary ergodic processes, the following corollary 
holds true.
\begin{corollary} There is no asymptotically correct discrimination procedure for the class of 
  stationary ergodic processes.
\end{corollary}

{\em Proof of Theorem~\ref{th}:}
We will assume that asymptotically correct discrimination  procedure $D$ for the class of all $B$-processes exists, and will construct a $B$-process $\rho$ such that 
if both sequences $x_i$ and $y_i$, $i\in\N$ are generated by $\rho$ then $\E D_n$ diverges; this contradiction will prove the theorem.

The scheme of the proof is as follows. 
On Step 1 we construct  a sequence of processes $\rho_{2k}$,  $\rho_{d2k+1}$, and $\rho_{u 2k+1}$, where $k=0,1,\dots$.
On Step 2 we construct a process $\rho$, which is shown to be  the limit of  the sequence $\rho_{2k}$, $k\in\N$, in $\bar d$-distance. 
On Step 3 we show  that two independent runs of the process $\rho$ have 
a property that (with high probability) they first behave like two runs  of a single process $\rho_0$, then like two runs
of two different processes $\rho_{u1}$ and $\rho_{d1}$, then like two runs of a single process $\rho_2$, and so on, thereby showing
that the test $D$ diverges and obtaining the desired contradiction.


Assume that there exists an asymptotically correct discriminating procedure~$D$. 
Fix some $\epsilon\in(0,1/2)$ and $\delta \in[1/2,1)$, to be defined on Step~3.

{\em Step 1.} We will construct the sequence of process $\rho_{2k}$,  $\rho_{u 2k+1}$, and $\rho_{d 2k+1}$, where $k=0,1,\dots$.

{\em Step 1.0.}
Construct the process $\rho_0$ as follows. A Markov chain $m_0$ is defined on the set $\N$ of states.
From each state $i\in\N$ the chain passes to the state $0$ with probability $\delta$ and to the state ${i+1}$ with
probability $1-\delta$. 
 With transition probabilities so defined, the chain possesses a unique stationary distribution $M_0$ on the set $\N$, 
 which can be calculated explicitly using e.g. \cite[Theorem VIII.4.1]{shir}, and is as follows:  $M_0(0)=\delta$,  $M_0(k)=\delta(1-\delta)^k$, 
  for all  $k\in\N$.
Take this distribution as the initial distribution over the states. 

The function $f_0$ maps the states to the output alphabet $\{0,1\}$ as follows: $f_0(i)=1$ for every $i\in\N$.
 Let $s_t$ be the state of the chain at time $t$. The process $\rho_0$
is defined as $\rho_0= f_0(s_t)$ for  $t\in\N$. 
As a result
of this definition, the process $\rho_0$ simply outputs $1$ with probability $1$ on every time step (however, by using different functions $f$ 
we will have less trivial processes in the sequel).  Clearly, the constructed process is stationary ergodic and a B-process.
So, we have defined the chain $m_0$ (and  the process $\rho_0$) up to a parameter $\delta$.


{\em Step 1.1.} We begin with the process $\rho_0$ and the chain $m_0$  of the previous step.
 Since the test D is asymptotically correct  we will have
$$
 \E_{\rho_0\times\rho_0}  D_{t_0}((x_1,\dots,x_{t_0}), (y_1,\dots,y_{t_0})) <\epsilon,
$$ from some  $t_0$ on, where both samples
$x_i$ and $y_i$ are generated by $\rho_0$ (that is, both samples consist of 1s only). 
Let $k_0$ be such an index that 
the chain $m_0$ starting from the state $0$ with probability $1$ does not reach the state $k_0-1$ by time $t_0$ (we can take $k_0=t_0+2$). 

Construct two processes $\rho_{u1}$ and $\rho_{d1}$ as follows. 
They are also based on the Markov chain $m_0$, but the functions $f$ are different. 
The function $f_{u1}:\N\rightarrow\{0,1\}$ is defined as follows: $f_{u1}(i)=f_0(i)=1$ for $i\le k_0$ and  $f_{u1}(i)=0$ for $i>k_0$.
The function $f_{d1}$ is identically $1$ ($f_{d1}(i)=1$, $i\in\N$).
The processes $\rho_{u1}$ and $\rho_{d1}$
are defined as $\rho_{u1}= f_{u1}(s_t)$  and $\rho_{d1}= f_{d1}(s_t)$ for $t\in\N$. Thus the process $\rho_{d1}$ will again produce only 1s, but the process $\rho_{u1}$ will occasionally produce 0s.

{\em Step 1.2.}
Being run on two samples generated by the processes $\rho_{u1}$ and $\rho_{d1}$ which both start from the state 0, 
 the test $D_n$ on the first $t_0$ steps produces many 0s,
since on these first $k_0$ states all the functions $f$, $f_{u1}$ and $f_{d1}$ coincide.
 However, since the processes are different and the test is asymptotically correct (by assumption), the test starts producing 1s, until by  a certain 
 time step $t_1$ almost all answers are 1s. 
Next we will construct the process $\rho_2$ by ``gluing'' together  $\rho_{u1}$ and $\rho_{d1}$ and continuing them in such a way that, being 
run on two samples produced by $\rho_2$ the test first produces 0s (as if the samples were drawn from $\rho_0$), then, with probability close to 1/2 it will
produce many 1s (as if the samples were from $\rho_{u1}$ and $\rho_{d1}$) and then again 0s. 

The process $\rho_2$ is the pivotal point of the construction, so we give it in some detail.  On step 1.2a we present 
the construction of the process, and on step 1.2b we show that this process is a $B$-process by demonstrating that 
it  is equivalent to a (deterministic) function of a Markov chain.

{\em Step 1.2a.}
Let $t_1>t_0$ be such a time index that 
$$
 \E_{\rho_{u1}\times\rho_{d1}}  D_k((x_1,\dots,x_{t_1}), (y_1,\dots,y_{t_1})) >1-\epsilon,
$$ where the samples
$x_i$ and $y_i$ are generated by $\rho_{u1}$ and $\rho_{d1}$ correspondingly (the samples are generated independently; that is, the process are based on two 
independent copies of the Markov chain $m_0$).
 Let $k_1>k_0$ be such an index that 
the chain $m$ starting from the state 0 with probability $1$ does not reach the state $k_1-1$ by time $t_1$.

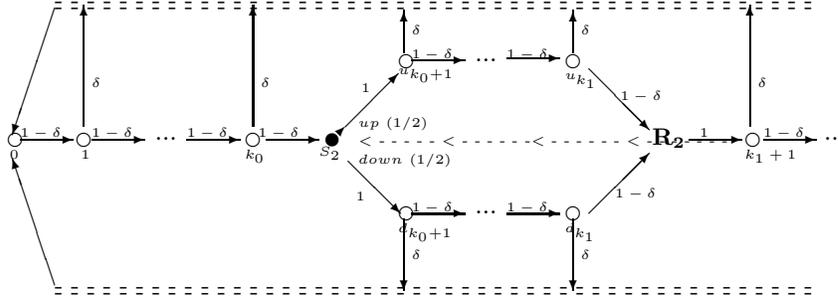
\begin{figure*}[h]\caption{The processes $m_2$ and $\rho_2$. {\small The states are depicted as circles, the arrows symbolize transition probabilities: from every state the process returns to 0 with probability $\delta$ or goes to the next state with probability $1-\delta$. From the switch $S_2$ the process
passes to the state indicated by the switch (with probability 1); here it is the state  $u_{k_0+1}$. When the process passes through the reset { $\bf R_2$} the switch $S_2$ is 
set to either $up$ or $down$ with equal probabilities. (Here $S_2$ is in the position $up$.) The function $f_2$ is 1 on all states except $u_{k0+1},\dots,u_{k1}$ where it is 0; $f_2$ applied to the states output by $m_2$ defines $\rho_2$.}}\label{fig:rho2}
\begin{picture}(200,120)(-70,-50)
\put(20,10){\circle{5}}
\put(18,3){\text{\tiny{0}}}
\put(22,10){\vector(1,0){20}}
\put(22,11){\text{\tiny{$1-\delta$}}}

\put(46,10){\circle{5}}
\put(45,3){\text{\tiny{1}}}
\put(49,10){\vector(1,0){20}}
\put(49,11){\text{\tiny{$1-\delta$}}}

\put(73,10){\text{{...}}}
\put(85,10){\vector(1,0){20}}
\put(85,11){\text{\tiny{$1-\delta$}}}
\put(110,10){\circle{5}}
\put(107,3){\text{\tiny{$k_0$}}}

\put(115,10){\vector(1,0){20}}
\put(112,11){\text{\tiny{$1-\delta$}}}
\put(135,4){\text{\tiny{$S_2$}}}
\put(140,10){\circle*{5}}
\put(143,13){\vector(1,1){2}}

\put(150,15){\text{\tiny{$up\ (1/2)$}}}
\put(150,1){\text{\tiny{$down\ (1/2)$}}}

\put(145,15){\vector(1,1){20}}
\put(151,28){\text{\tiny{1}}}

\put(168,40){\circle{5}}
\put(165,35){\text{\tiny{$u_{k_0+1}$}}}
\put(170,40){\vector(1,0){20}}
\put(170,41){\text{\tiny{$1-\delta$}}}
\put(194,40){\text{{...}}}
\put(206,41){\vector(1,0){20}}
\put(206,41){\text{\tiny{$1-\delta$}}}
\put(231,40){\circle{5}}
\put(228,33){\text{\tiny{$u_{k_1}$}}}

\put(146, 2){\vector(1,-1){20}}
\put(149,-13){\text{\tiny{1}}}

\put(168,-18){\circle{5}}
\put(165,-25){\text{\tiny{$d_{k_0+1}$}}}
\put(170,-18){\vector(1,0){20}}
\put(170,-17){\text{\tiny{$1-\delta$}}}
\put(194,-18){\text{{...}}}
\put(206,-18){\vector(1,0){20}}
\put(206,-17){\text{\tiny{$1-\delta$}}}
\put(231,-18){\circle{5}}
\put(228,-25){\text{\tiny{$d_{k_1}$}}}

\put(237,37){\vector(1,-1){23}}
\put(249,25){\text{\tiny{$1-\delta$}}}
\put(247,-12){\text{\tiny{$1-\delta$}}}
\put(237,-18){\vector(1,1){23}}

\put(261,8){\text{\small{$\mathbf{R_2}$}}}
\put(150,8){\text{\tiny{$<$ - - - - - $<$ - - - - - -$<$ - - - - - - $<$ - - - - - - -}}}

\put(275,10){\vector(1,0){20}}
\put(279,11){\text{\tiny{$1$}}}
\put(299,10){\circle{5}}
\put(296,3){\text{\tiny{$k_1+1$}}}
\put(303,10){\vector(1,0){20}}
\put(303,11){\text{\tiny{$1-\delta$}}}
\put(326,10){\text{{...}}}

\put(35,60){\text{\small{- - - - - - - - - - - - - - - - - - - - - - - - - - - - - - - - - - - - - - - - - - - - - - -}}}
\put(35,58){\text{\small{- - - - - - - - - - - - - - - - - - - - - - - - - - - - - - - - - - - - - - - - - - - - - - -}}}
\put(46,15){\vector(0,1){46}}
\put(49,30){\text{\tiny{$\delta$}}}
\put(110,15){\vector(0,1){46}}
\put(113,30){\text{\tiny{$\delta$}}}

\put(167,43){\vector(0,1){16}}
\put(170,50){\text{\tiny{$\delta$}}}
\put(231,43){\vector(0,1){16}}
\put(234,50){\text{\tiny{$\delta$}}}

\put(167,-20){\vector(0,-1){27}}
\put(170,-35){\text{\tiny{$\delta$}}}
\put(231,-20){\vector(0,-1){27}}
\put(234,-35){\text{\tiny{$\delta$}}}

\put(298,15){\vector(0,1){46}}
\put(301,30){\text{\tiny{$\delta$}}}

\put(35,60){\vector(-1,-3){16}}
\put(35,-45){\vector(-1,3){16}}

\put(35,-50){\text{\small{- - - - - - - - - - - - - - - - - - - - - - - - - - - - - - - - - - - - - - - - - - - - - - -}}}
\put(35,-48){\text{\small{- - - - - - - - - - - - - - - - - - - - - - - - - - - - - - - - - - - - - - - - - - - - - - -}}}
\end{picture}
\end{figure*}

Construct the process $\rho_2$ as follows (see fig.~\ref{fig:rho2}).
 It is based on a chain $m_2$ on which Markov assumption is violated. 
The transition probabilities on states $0,\dots,k_0$ are the same as for the Markov chain $m$ (from each
state return to 0 with probability $\delta$ or go to the next state with probability $1-\delta$). 

There are two ``special'' states: the ``switch'' $S_2$ and 
the ``reset'' $R_2$.  From the state $k_0$ the chain passes with probability $1-\delta$  to the ``switch'' state
$S_2$. 
The switch $S_2$ can itself have two values: $up$ and $down$. If $S_2$ has the value $up$ then from $S_2$ the chain passes to the state $u_{k_0+1}$ with probability 1, while
if $S_2=down$ the chain goes to $d_{k_0+1}$, with probability 1.  If the chain reaches the state $R_2$ then 
the value of $S_2$ is set to $up$ with probability 1/2 and with probability 1/2 it is set to $down$.
 In other words, the first transition from $S_2$ is random 
(either to $u_{k_0+1}$ or to  $d_{k_0+1}$ with equal probabilities) and then this decision is remembered until the ``reset'' state $R_2$ is visited, whereupon the
switch again assumes the values $up$ and $down$ with equal probabilities.

The rest of the transitions are as follows. From each state $u_i$, $k_0\le i\le k_1$ the chain passes to the state $0$ with probability $\delta$ and to the next state $u_{i+1}$ with probability $1-\delta$.
From the state $u_{k_1}$ the process goes with probability $\delta$ to 0 and with probability $1-\delta$ to the ``reset'' state $R_2$. The same with states $d_i$: 
for $k_0<i\le k_1$ the process returns to 0 with probability $\delta$ or goes to the next state $d_{i+1}$ with probability $1-\delta$, where the next state
for $d_{k_1}$ is the ``reset'' state $R_2$. From $R_2$ the process goes with probability 1 to the state $k_1+1$ where from the chain continues ad infinitum:
to the state 0 with probability $\delta$ or to the next state $k_1+2$ etc. with probability $1-\delta$.

The initial distribution on the states is defined as follows.
 The probabilities of the states $0..k_0, k_1+1,k_1+2,\dots$ are the same
as in the Markov chain $m_0$, that is, $\delta(1-\delta)^j$, for $j=0..k_0, k_1+1,k_1+2,\dots$. 
For the states $u_j$ and $d_j$, $k_0<j\le k_1$ define  their initial probabilities to be 1/2 of the probability of the corresponding
state in the chain $m_0$, that is $m_2(u_j)=m_2(d_j)=m_0(j)/2=\delta(1-\delta)^j/2$. Furthermore, if the chain starts in a state $u_j$,
$k_0<j\le k_1$, then the value of the switch $S_2$ is $up$, and if it starts in the state $d_j$ then the value of the switch $S_2$ is $down$,
whereas if the chain starts in any other state then the probability distribution on the values of the switch $S_2$ is  1/2 for either $up$ or $down$.

The function $f_2$ is defined as follows: $f_2(i)=1$ for $0\le i\le k_0$ and $i>k_1$ (before the switch and after the reset); $f_2(u_i)=0$ for all $i$, $k_0<i\le k_1$ and $f_2(d_i)=1$ for all $i$, $k_0<i\le k_1$.
The function $f_2$ is undefined on $S_2$ and $R_2$, therefore there is no output on these states (we also assume that passing through $S_2$ and $R_2$ does not increment time). As before, the process $\rho_2$ is defined as $\rho_2=f_2(s_t)$ where
$s_t$ is the state of $m_2$ at time $t$, omitting the states $S_2$ and $R_2$.
The resulting process s illustrated on fig.~\ref{fig:rho2}.

\begin{figure*}[!t]\caption{The process $m_2'$. \small{The function $f_2$ is 1 everywhere except the states $u_{k_0+1},\dots,u_{k_1}$, where it is 0.}}\label{fig:m2}
\begin{picture}(200,120)(-70,-50)




\put(20,40){\circle{5}}
\put(18,34){\text{\tiny{$u_0$}}}
\put(22,40){\vector(1,0){20}}
\put(22,41){\text{\tiny{$1-\delta$}}}

\put(46,40){\circle{5}}
\put(45,34){\text{\tiny{$u_1$}}}
\put(49,40){\vector(1,0){20}}
\put(49,41){\text{\tiny{$1-\delta$}}}

\put(73,40){\text{{...}}}
\put(85,40){\vector(1,0){20}}
\put(85,41){\text{\tiny{$1-\delta$}}}
\put(110,40){\circle{5}}
\put(107,34){\text{\tiny{$u_{k_0-1}$}}}

\put(115,40){\vector(1,0){20}}
\put(112,41){\text{\tiny{$1-\delta$}}}
\put(140,40){\circle{5}}
\put(137,34){\text{\tiny{$u_{k_0}$}}}
\put(143,40){\vector(1,0){20}}
\put(142,41){\text{\tiny{$1-\delta$}}}

\put(20,-20){\circle{5}}
\put(18,-26){\text{\tiny{$d_0$}}}
\put(22,-20){\vector(1,0){20}}
\put(22,-19){\text{\tiny{$1-\delta$}}}

\put(46,-20){\circle{5}}
\put(45,-26){\text{\tiny{$d_1$}}}
\put(49,-20){\vector(1,0){20}}
\put(49,-19){\text{\tiny{$1-\delta$}}}

\put(73,-20){\text{{...}}}
\put(85,-20){\vector(1,0){20}}
\put(85,-19){\text{\tiny{$1-\delta$}}}
\put(110,-20){\circle{5}}
\put(107,-26){\text{\tiny{$d_{k_0-1}$}}}

\put(115,-20){\vector(1,0){20}}
\put(112,-19){\text{\tiny{$1-\delta$}}}
\put(140,-20){\circle{5}}
\put(137,-26){\text{\tiny{$d_{k_0}$}}}
\put(143,-20){\vector(1,0){20}}
\put(142,-19){\text{\tiny{$1-\delta$}}}


\put(168,40){\circle{5}}
\put(165,35){\text{\tiny{$u_{k_0+1}$}}}
\put(170,40){\vector(1,0){20}}
\put(170,41){\text{\tiny{$1-\delta$}}}
\put(194,40){\text{{...}}}
\put(206,41){\vector(1,0){20}}
\put(206,41){\text{\tiny{$1-\delta$}}}
\put(231,40){\circle{5}}
\put(228,33){\text{\tiny{$u_{k_1}$}}}


\put(168,-18){\circle{5}}
\put(165,-25){\text{\tiny{$d_{k_0+1}$}}}
\put(170,-18){\vector(1,0){20}}
\put(170,-17){\text{\tiny{$1-\delta$}}}
\put(194,-18){\text{{...}}}
\put(206,-18){\vector(1,0){20}}
\put(206,-17){\text{\tiny{$1-\delta$}}}
\put(231,-18){\circle{5}}
\put(228,-25){\text{\tiny{$d_{k_1}$}}}

\put(237,37){\vector(1,-1){23}}
\put(249,25){\text{\tiny{$1-\delta$}}}
\put(247,-12){\text{\tiny{$1-\delta$}}}
\put(237,-18){\vector(1,1){23}}

\put(265,10){\circle{5}}
\put(262,3){\text{\tiny{$k_1+1$}}}

\put(270,10){\vector(1,0){20}}
\put(270,11){\text{\tiny{$1-\delta$}}}
\put(294,10){\circle{5}}
\put(291,3){\text{\tiny{$k_1+2$}}}
\put(298,10){\vector(1,0){20}}
\put(298,11){\text{\tiny{$1-\delta$}}}
\put(321,10){\text{{...}}}

\put(35,60){\text{\small{- - - - - - - - - - - - - - - - - - - - - - - - - - - - - - - - - - - - - - - - - - - - - - }}}
\put(35,58){\text{\small{- - - - - - - - - - - - - - - - - - - - - - - - - - - - - - - - - - - - - - - - - - - - - - }}}
\put(46,45){\vector(0,1){16}}
\put(48,49){\text{\tiny{$\delta$}}}
\put(110,45){\vector(0,1){16}}
\put(112,49){\text{\tiny{$\delta$}}}

\put(167,43){\vector(0,1){16}}
\put(170,50){\text{\tiny{$\delta$}}}
\put(231,43){\vector(0,1){16}}
\put(234,50){\text{\tiny{$\delta$}}}

\put(167,-20){\vector(0,-1){27}}
\put(170,-35){\text{\tiny{$\delta$}}}
\put(231,-20){\vector(0,-1){27}}
\put(234,-35){\text{\tiny{$\delta$}}}

\put(109,-26){\vector(0,-1){21}}
\put(110,-35){\text{\tiny{$\delta$}}}
\put(139,-26){\vector(0,-1){21}}
\put(141,-35){\text{\tiny{$\delta$}}}
\put(45,-26){\vector(0,-1){21}}
\put(47,-35){\text{\tiny{$\delta$}}}

\put(35,-50){\text{\small{- - - - - - - - - - - - - - - - - - - - - - - - - - - - - - - - - - - - - - - - - - - - - - }}}
\put(35,-48){\text{\small{- - - - - - - - - - - - - - - - - - - - - - - - - - - - - - - - - - - - - - - - - - - - - - }}}

\put(294,13){\vector(0,1){45}}
\put(295,34){\text{\tiny{$\delta/2$}}}
\put(265,13){\vector(0,1){45}}
\put(266,34){\text{\tiny{$\delta/2$}}}

\put(294,3){\vector(0,-1){48}}
\put(295,-34){\text{\tiny{$\delta/2$}}}
\put(265,3){\vector(0,-1){48}}
\put(266,-34){\text{\tiny{$\delta/2$}}}



\put(35,60){\vector(-1,-1){16}}
\put(38,-45){\vector(-1,1){19}}
\end{picture}
\end{figure*}
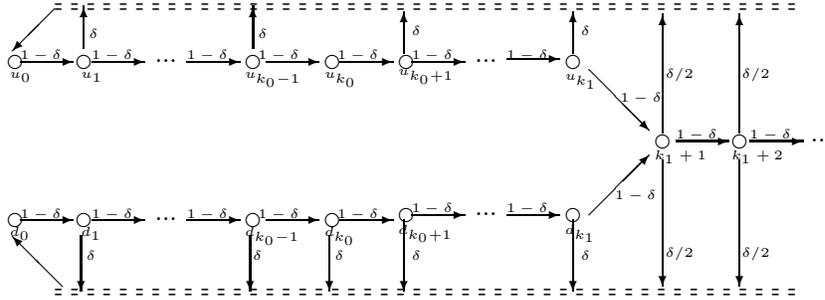

{\em Step 1.2b.}
To show that the process $\rho_2$ is stationary ergodic and a $B$-process, we will 
show that it is equivalent to a function of a stationary ergodic Markov chain, whereas all such process
are known to be $B$ (e.g. \cite{shbook}). 
The construction is as follows  (see fig.~\ref{fig:m2}). This chain has states 
 $k_1+1,\dots$ and also 
$u_0,\dots,u_{k_0},u_{k_0+1},\dots,u_{k_1}$ and
$d_0,\dots,d_{k_0},d_{k_0+1},\dots,d_{k_1}$.
 From the states $u_i$, $i=0,\dots,k_1$ the 
chain passes with probability $1-\delta$ to the next state $u_{i+1}$, where the next state for $u_{k_1}$ is $k+1$ and with probability $\delta$ returns to the state $u_0$ (and not to the state 0). Transitions
for the state $d_0,\dots,d_{k_1-1}$ are defined analogously. Thus the states $u_{k_i}$ correspond to the state $up$ of the switch $S_2$ and the states $d_{k_i}$~--- to the state
$down$ of the switch. Transitions for the states $k+1,k+2,\dots$ are defined as follows: with probability $\delta/2$ to the state $u_0$, with probability $\delta/2$ to 
the state $d_0$,  and with probability $1-\delta$ to the next state.
Thus, transitions to 0 from  the states with indices greater than $k_1$ corresponds to the reset $R_2$. Clearly, the chain $m_2'$ as defined possesses a unique stationary distribution $M_2$ over the set of states and $M_2(i)>0$ for every state $i$. Moreover, this distribution is the same as the initial distribution on the states of the chain $m_0$,
except  for the states $u_i$ and $d_i$, for which we have $m_2'(u_i)=m_2'(d_i)=m_0(i)/2=\delta(1-\delta)^i/2$, for $0\le i\le k_0$. 
We take  this distribution as its initial distribution on the states of $m_2'$. The resulting process $m_2'$ is stationary ergodic, and a $B$-process, since
it is   a function of a Markov chain \cite{shbook}. 
It is easy to see that if we define the function $f_2$ on the states of $m_2'$ as 1 on all states except $u_{k_0+1},\dots,u_{k_1}$, then the resulting
process is exactly the process $\rho_2$. Therefore, $\rho_2$ is stationary ergodic and a $B$-process.

{\em Step 1.$k$.}
As before, we can continue the construction of 
 the processes $\rho_{u3}$ and $\rho_{d3}$, that start with a segment of $\rho_2$. 
Let $t_2>t_1$ be a time index such that 
$$
 \E_{\rho_2\times\rho_2}  D_{t_2} <\epsilon,
$$ where both samples are generated by $\rho_2$. Let $k_2>k_1$ be  such an index that when starting from the state 0  the process $m_2$ with probability 1 does not 
reach $k_2-1$ by time $t_2$ (equivalently: the process $m_2'$ does not reach $k_2-1$ when starting from either 0, $u_0$ or $d_0$). 
The processes $\rho_{u3}$ and $\rho_{d3}$ are based on the same process $m_2$ as  $\rho_2$. The functions $f_{u3}$ and $f_{d3}$ coincide with $f_2$  on all states up to the state $k_2$ (including the states $u_i$ and $d_i$, $k_0<i\le k_1$).  After $k_2$ the function $f_{u3}$ outputs 0s while $f_{d3}$ outputs 1s: $f_{u3}(i)=0$, $f_{d3}(i)=1$ for 	$i>k_2$.

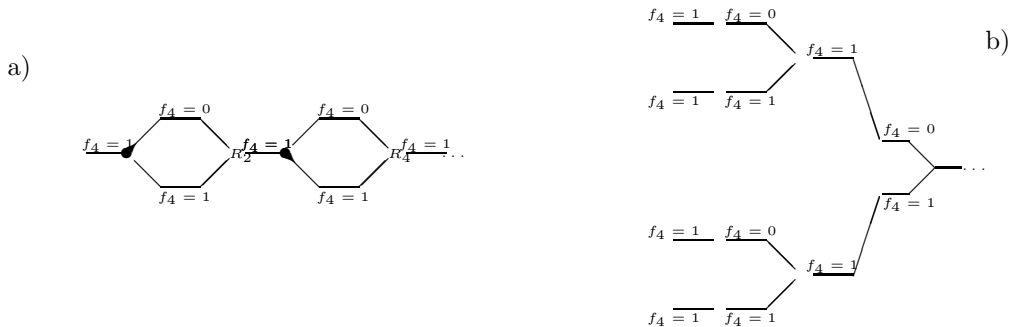
\begin{figure*}[h]\caption{a) The processes $m_4$. 
b) The Markov chain $m_4'$}\label{fig:m4}
\begin{picture}(200,120)(0,-50)

\newsavebox{\blo}
\savebox{\blo}(0,0)[l]{
\put(20,10){\line(1,0){15}}
\put(35,10){\circle*{4}}
\put(37,12){\vector(1,1){2}}

\put(38,13){\line(1,1){10}}
\put(38,7){\line(1,-1){10}}

\put(48,23){\line(1,0){15}}
\put(48,-3){\line(1,0){15}}

\put(63,23){\line(1,-1){11}}
\put(63,-3){\line(1,1){11}}
\put(74,8){\text{\tiny{$R_2$}}}

\put(81,10){\line(1,0){15}}
\put(47,25){\text{\tiny{$f_4=0$}}}
\put(18,12){\text{\tiny{$f_4=1$}}}
\put(47,-8){\text{\tiny{$f_4=1$}}}
\put(78,12){\text{\tiny{$f_4=1$}}}

}

\newsavebox{\blc}
\savebox{\blc}(0,0)[l]{
\put(28,23){\line(1,0){15}}
\put(28,-3){\line(1,0){15}}

\put(48,23){\line(1,0){15}}
\put(48,-3){\line(1,0){15}}

\put(63,23){\line(1,-1){11}}
\put(63,-3){\line(1,1){11}}

\put(81,10){\line(1,0){15}}
\put(47,25){\text{\tiny{$f_4=0$}}}
\put(18,25){\text{\tiny{$f_4=1$}}}
\put(18,-8){\text{\tiny{$f_4=1$}}}
\put(47,-8){\text{\tiny{$f_4=1$}}}
\put(78,12){\text{\tiny{$f_4=1$}}}

}

\newsavebox{\bloo}
\savebox{\bloo}(0,0)[l]{
\put(20,10){\line(1,0){15}}
\put(35,10){\circle*{4}}
\put(37,8){\vector(1,-1){2}}

\put(38,13){\line(1,1){10}}
\put(38,7){\line(1,-1){10}}

\put(48,23){\line(1,0){15}}
\put(48,-3){\line(1,0){15}}

\put(63,23){\line(1,-1){11}}
\put(63,-3){\line(1,1){11}}
\put(74,8){\text{\tiny{$R_4$}}}

\put(81,10){\line(1,0){15}}
\put(47,25){\text{\tiny{$f_4=0$}}}
\put(18,12){\text{\tiny{$f_4=1$}}}
\put(47,-8){\text{\tiny{$f_4=1$}}}
\put(78,12){\text{\tiny{$f_4=1$}}}

}

\put(40,10){\usebox{\blo}}

\put(100,10){\usebox{\bloo}}

\put(194,9){\text{\tiny{$\dots$}}}
%






\put(254,46){\usebox{\blc}}
\put(254,-36){\usebox{\blc}}%

\put(361,15){\line(1,0){10}}
\put(361,-5){\line(1,0){10}}
\put(371,-5){\line(1,1){10}}
\put(371,15){\line(1,-1){10}}
\put(381,5){\line(1,0){10}}
\put(391,4){\text{\tiny{$\dots$}}}

\put(350,46){\line(1,-3){10}}
\put(350,-36){\line(1,3){10}}

\put(361,18){\text{\tiny{$f_4=0$}}}
\put(361,-10){\text{\tiny{$f_4=1$}}}

\put(30,40){\text{a)}}


\put(400,50){\text{b)}}

\end{picture}
\end{figure*}

Furthermore, we find a time  $t_3>t_2$ by which we have 
$
 \E_{\rho_{u3}\times\rho_{d3}}  D_{t_3} >1-\epsilon,
$ where the samples are generated by  $\rho_{u3}$ and $\rho_{d3}$, which is possible since $D$ is consistent.
Next, find an  index $k_3>k_2$ such that the process $m_2$ does not reach $k_3-1$ with probability $1$ if the processes $\rho_{u3}$ and $\rho_{d3}$ 
are used to produce two independent sequences and both start from the state 0.
  We then construct the process $\rho_4$ based on a (non-Markovian) process $m_4$ by ``gluing''
together $\rho_{u3}$ and $\rho_{d3}$ after the step $k_3$ with a switch $S_4$ and a reset $R_4$ exactly as was done when constructing the process $\rho_2$. 
The process $m_4$ is illustrated on fig.~\ref{fig:m4}a). 
The process $m_4$ can be shown to be equivalent to a Markov chain $m_4'$, 
 which  is constructed analogously to the chain $m_2'$  
(see fig.~\ref{fig:m4}b).
Thus, the process $\rho_4$ is can be shown to be a $B$-process.

Proceeding this way we can construct the processes $\rho_{2j}$, $\rho_{u2j+1}$ and $\rho_{d2j+1}$, $j\in\N$ choosing the time steps $t_j>t_{j-1}$ so that
the expected output of the test approaches 0 by the time $t_j$ being run on  two samples produced by $\rho_j$ for even $j$, and approaches 1 by the time $t_j$ being run on samples produced by $\rho_{uj}$ and $\rho_{dj}$ for odd $j$:
\begin{equation}\label{eq:teven}
 \E_{\rho_{2j}\times\rho_{2j}}  D_{t_{2j}} <\epsilon
\end{equation}
and 
\begin{equation}\label{eq:todd}
 \E_{\rho_{u2j+1}\times\rho_{d2j+1}}  D_{t_{2j+1}} > (1-\epsilon).
\end{equation}
For each $j$ the number $k_j>k_{j-1}$  is selected in a such a way that the state  $k_j-1$ is not reached (with probability 1) by the  time $t_j$ when starting from the state 0.
Each of the processes $\rho_{2j}$, $\rho_{u2j+1}$ and $\rho_{dj2+1}$, $j\in\N$ can be shown to be stationary ergodic and a $B$-process by demonstrating
equivalence to a Markov chain, analogously to the Step 1.2. The initial state distribution of each of the processes $\rho_t, t\in\N$ is  $M_{t}(k)=\delta(1-\delta)^k$ and 
$M_{t}(u_k)=M_{t}(d_k)= \delta(1-\delta)^k/2$ for those $k\in\N$ for which the corresponding states are defined.

{\em Step 2.} 
Having defined $k_j$, $j\in\N$ we can define the process $\rho$. The construction is given on Step~2a, while on Step~2b we show
that $\rho$ is stationary ergodic and a $B$-process, by showing that it is the limit of the sequence  $\rho_{2j}$, $j\in\N$. 

\begin{figure}[h]\caption{The processes $m_\rho$ and $\rho$. The states are on horizontal lines. The function $f$ being applied to the states of $m_\rho$ 
defines the process $\rho$. Its value is  $0$ on the states on the upper lines (states $u_{k_{2j}+1},\dots,u_{k_{2j+1}}$, where $k\in\N$) and 1 on the rest of the states.}\label{fig:rho}
\begin{picture}(200,47)(-140,-7)

\newsavebox{\bloc}
\savebox{\bloc}(0,0)[l]{
\put(20,10){\line(1,0){15}}
\put(35,10){\circle*{4}}
\put(37,12){\vector(1,1){2}}

\put(38,13){\line(1,1){10}}
\put(38,7){\line(1,-1){10}}

\put(48,23){\line(1,0){15}}
\put(48,-3){\line(1,0){15}}

\put(63,23){\line(1,-1){11}}
\put(63,-3){\line(1,1){11}}
\put(74,8){\text{\tiny{$R_2$}}}

\put(81,10){\line(1,0){15}}
\put(47,25){\text{\tiny{$f=0$}}}
\put(18,12){\text{\tiny{$f=1$}}}
\put(47,-8){\text{\tiny{$f=1$}}}

}
\newsavebox{\blocc}
\savebox{\blocc}(0,0)[l]{
\put(20,10){\line(1,0){15}}
\put(35,10){\circle*{4}}
\put(37,8){\vector(1,-1){2}}

\put(38,13){\line(1,1){10}}
\put(38,7){\line(1,-1){10}}

\put(48,23){\line(1,0){15}}
\put(48,-3){\line(1,0){15}}

\put(63,23){\line(1,-1){11}}
\put(63,-3){\line(1,1){11}}
\put(74,8){\text{\tiny{$R_4$}}}

\put(81,10){\line(1,0){15}}
\put(47,25){\text{\tiny{$f=0$}}}
\put(18,12){\text{\tiny{$f=1$}}}
\put(47,-8){\text{\tiny{$f=1$}}}

}

\newsavebox{\bloccc}
\savebox{\bloccc}(0,0)[l]{
\put(20,10){\line(1,0){15}}
\put(35,10){\circle*{4}}
\put(37,12){\vector(1,1){2}}

\put(38,13){\line(1,1){10}}
\put(38,7){\line(1,-1){10}}

\put(48,23){\line(1,0){15}}
\put(48,-3){\line(1,0){15}}

\put(63,23){\line(1,-1){11}}
\put(63,-3){\line(1,1){11}}
\put(74,8){\text{\tiny{$R_6$}}}

\put(81,10){\line(1,0){15}}
\put(47,25){\text{\tiny{$f=0$}}}
\put(18,12){\text{\tiny{$f=1$}}}
\put(47,-8){\text{\tiny{$f=1$}}}

}

\put(-10,10){\usebox{\bloc}}
\put(50,10){\usebox{\blocc}}
\put(110,10){\usebox{\bloccc}}
\put(212,10){\text{\tiny{$\dots$}}}

\end{picture}
\end{figure}
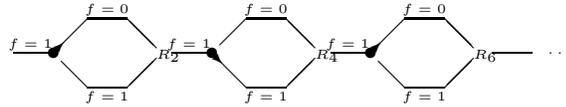

{\em Step 2a.}
The process $\rho$ can be constructed as follows (see fig.~\ref{fig:rho}). 
 The construction is based on the (non-Markovian) process $m_\rho$ that has states $0,\dots,k_0$, $k_{2j+1}+1,\dots,k_{2(j+1)}$, $u_{k_{2j}+1},\dots,u_{k_{2j+1}}$ and $d_{k_{2j}+1},\dots,d_{k_{2j+1}}$
for $j\in\N$, along with switch states $S_{2j}$ and reset states $R_{2j}$. Each switch $S_{2j}$ diverts the process to the state $u_{k_{2j}+1}$ if the switch has value $up$ and
to $d_{k_{2j}+1}$ if it has the value $down$. The reset $R_{2j}$ sets $S_{2j}$ to  $up$ with probability 1/2 and to $down$ also with probability 1/2.
 From each state that is neither a reset nor a switch, the process goes to the next state with probability $1-\delta$ and returns
to the state 0 with probability $\delta$ (cf. Step 1$k$). 

The initial distribution $M_\rho$ on the states of $m_\rho$ is defined as follows. For every state $i$ such that $0 \le i\le k_0$ and  $k_{2j+1}<i\le k_{2_j+2}$, $j=0,1,\dots$,
define the initial probability of the state $i$ as $M_{\rho}(i)=\delta(1-\delta)^i$ (the same as in the chain $m_0$), and for the sets $u_j$ and $d_j$ 
(for those $j$ for which these sets are defined)
let $M_{\rho}(u_j)=M_{\rho}(d_j):=\delta(1-\delta)^i/2$ (that is, 1/2 of the probability of the corresponding state of $m_0$). 

The function $f$ is defined as 1 everywhere except 
for the states $u_j$ (for all $j\in\N$ for which $u_j$ is defined)  on which $f$ takes the value 0.  The process $\rho$ is defined at time $t$ as $f(s_t)$,
where $s_t$ is the state of $m_\rho$ at time $t$. 

{\em Step 2b.} To show that $\rho$ is a $B$-process, let us first show that it is stationary. 
To do this, define the so-called distributional distance on the set of all stochastic processes as follows. 
$$
d(\mu_1,\mu_2)=\sum_{i=1}^\infty w_i |\mu_1((x_{1},\dots, x_{|B_i|})=B_i)-\mu_2((x_{1},\dots, x_{|B_i|})=B_i)|,
$$
where $\mu_1,\mu_2$ are any stochastic processes,  $w_k:=2^{-k}$ and $B_i$ ranges over all tuples $B\in \cup_{k\in\N}X^k$, assuming
some fixed order on this set. The set of all stochastic processes, equipped with this distance, is complete, and the set of all
stationary processes is its closed subset \cite{gray}. 
Thus, to show that the process $\rho$ is stationary it suffices to show that $\lim_{j\to\infty}d(\rho_{2j},\rho)=0$, since
the processes $\rho_{2j}$, $j\in\N$, are stationary. To do this, it is enough to demonstrate
that 
\begin{equation}\label{eq:lim}
\lim_{j\to\infty}  |\rho((x_{1},\dots, x_{|B|})=B)-\rho_{2j}((x_{1},\dots, x_{|B|})=B)|=0
\end{equation}
for each $B\in \cup_{k\in\N}X^k$.
  Since the processes $m_\rho$ and $m_{2j}$ coincide on all 
states up to $k_{2j+1}$, we have 
$$
|\rho(x_{n}=a)-\rho_{2j}(x_{n}=a)| = |\rho(x_{1}=a)-\rho_{2j}(x_{1}=a)|=\le \sum_{k>k_{2j+1}} M_\rho(k) + \sum_{k>k_{2j+1}} M_{2j}(k) 
$$ for every $n\in\N$ and $a\in X$. 
Moreover, 
 for any  tuple $B\in\cup_{k\in\N}X^k$ we obtain
$$
  |\rho((x_{1},\dots, x_{|B|})=B)-\rho_{2j}((x_{1},\dots, x_{|B|})=B)|\le  |B|\left(\sum_{k>k_{2j+1}} M_\rho(k) + \sum_{k>k_{2j+1}} M_{2j}(k)\right)\to0
$$
where the convergence follows from $k_{2j}\to\infty$. 
We conclude that~(\ref{eq:lim}) holds true, so that $d(\rho,\rho_{2j})\to 0$ and  $\rho$ is stationary.

To show that $\rho$  is a $B$-process, we will demonstrate that it
is the limit of the sequence $\rho_{2k}$, $k\in\N$ in the $\bar d$ distance (which was only defined for stationary processes). Since the set of all 
$B$-process is a closed subset of all stationary processes, it will follow that $\rho$ itself
is a $B$-process.  (Observe that this way we  get ergodicity of $\rho$ ``for free'', since the set of all ergodic processes is closed in $\bar d$ distance,
and all the processes $\rho_{2j}$ are ergodic.)
In order to show that $\bar d(\rho,\rho_{2k})\to0$ we have
 to find for each $j$ a  processes $\nu_{2j}$ on pairs $(x_1,y_1),(x_2,y_2),\dots$, such 
that $x_i$ are distributed according to $\rho$ and $y_i$ are distributed according to $\rho_{2j}$, 
and such that $\lim_{j\to\infty}\nu_{2j}(x_1\ne y_1)=0$. Construct such a coupling as follows.
Consider the chains $m_\rho$ and $m_{2j}$, which start in the same state (with initial distribution being $M_\rho$) and always take state transitions together,
where if the process $m_\rho$ is in the  state $u_{t}$ or $d_t$, $t\ge k_{2j+1}$ (that is, one of the states which the chain $m_{2j}$ does not have) 
then the chain $m_{2j}$ is in  the state $t$. The first coordinate of the process $\nu_{2j}$ is obtained by applying the function $f$ to the
process $m_\rho$ and the second by applying $f_{2j}$ to the chain $m_{2j}$. Clearly, the distribution 
of the first coordinate is $\rho$ and the distribution of the second is $\rho_{2j}$. Since the chains start in the same 
state and always take state transitions together, and since the chains $m_\rho$ and $m_{2j}$ coincide up to the 
state $k_{2j+1}$  we have $\nu_{2j}(x_1\ne y_1)\le \sum_{k>k_{2j+1}} M_\rho(k)\to0$. Thus, $\bar d(\rho, \rho_{2j})\to0$, so that  $\rho$ is a $B$-process.

{\em Step 3.}
Finally, it remains to show that the expected output of the test $D$ diverges if the test is run on two independent samples produced by $\rho$.

Recall that for all the chains $m_{2j}$, $m_{u2j+1}$ and $m_{d2j+1}$  as well as for the chain  $m_\rho$, the initial
probability of the state 0 is  $\delta$. By construction, if the process $m_\rho$ starts at the state 0 then up to the time step $k_{2j}$ it behaves exactly 
as $\rho_{2j}$ that has started at the state 0.  In symbols, we have 
\begin{equation}\label{eq:even1}
 E_{\rho\times\rho} ( D_{t_{2j}}| s_0^x=0, s_0^y=0)=E_{\rho_{2j}\times\rho_{2j}} ( D_{t_{2j}}| s_0^x=0, s_0^y=0)
\end{equation}
for $j\in\N$, where $s_0^x$ and $s_0^y$ denote the initial states of the processes generating the samples $x$ and $y$ correspondingly.

We will use the following simple decomposition 
\begin{equation}\label{eq:dec}
 \E( D_{t_j}) = \delta^2\E( D_{t_j}|s^x_0=0, s^y_0=0) 
+ (1-\delta^2)\E( D_{t_j}|s^x_0\ne0 \text{ or } s^y_0\ne0),
\end{equation}
(\ref{eq:even1}), and~(\ref{eq:teven}) we have 
\begin{multline}\label{eq:even2}
 \E_{\rho\times\rho} (  D_{t_{2j}}) \le \delta^2 \E_{\rho\times\rho} ( D_{t_{2j}} | s^x_0=0, s^y_0=0) +(1-\delta^2) \\ 
= \delta^2 \E_{\rho_{2j}\times\rho_{2j}} ( D_{t_{2j}} | s^x_0=0, s^y_0=0) +(1-\delta^2) \\
\le \E_{\rho_{2j}\times\rho_{2j}} +(1-\delta^2)  < \epsilon +(1-\delta^2).
\end{multline}

 For odd indices, if the process $\rho$ starts at the state 0 then (from the definition of $t_{2j+1}$) by the  time $t_{2j+1}$ it does not 
reach the reset $R_{2j}$;  therefore, in this case  the value of the switch $S_{2j}$ does not change up to the time $t_{2j+1}$.
 Since the definition of $m_\rho$ is symmetric with respect to the values $up$ and $down$ of each switch, the probability that two samples $x_1,\dots,x_{t_{2j+1}}$ and $y_1,\dots,y_{t_{2j+1}}$ generated  independently
by (two runs of) the process $\rho$ produced different values of the switch $S_{2j}$ when passing through it for  the first time is 1/2. In other words, with probability 1/2 two samples generated 
by  $\rho$ starting at the state 0 will look by the time $t_{2j+1}$ as two samples generated by $\rho_{u2j+1}$ and $\rho_{d2j+1}$ that has started at state 0. Thus
\begin{equation}\label{eq:odd1}
 E_{\rho\times\rho} ( D_{t_{2j+1}}| s_0^x=0, s_0^y=0)\\ \ge \frac{1}{2}E_{\rho_{u2j+1}\times\rho_{d2j+1}} ( D_{t_{2j+1}}| s_0^x=0, s_0^y=0)
\end{equation}
for $j\in\N$. Using this, (\ref{eq:dec}),  and~(\ref{eq:todd}) we obtain
\begin{multline}\label{eq:odd2}
 \E_{\rho\times\rho} (  D_{t_{2j+1}}) \ge  \delta^2 \E_{\rho\times\rho} ( D_{t_{2j+1}} | s^x_0=0, s^y_0=0)\\
\ge \frac{1}{2} \delta^2 \E_{\rho_{2j+1}\times\rho_{2j+1}} ( D_{t_{2j+1}} | s^x_0=0, s^y_0=0)\\
\ge  \frac{1}{2} \left( \E_{\rho_{2j+1}\times\rho_{2j+1}} (  D_{t_{2j+1}}) - (1-\delta^2)\right)> \frac{1}{2} (\delta^2-\epsilon).
\end{multline}


Taking $\delta$ large and $\epsilon$ small (e.g. $\delta=0.9$ and $\epsilon=0.1$), we can make the bound~(\ref{eq:even2}) close to 0 and the bound~(\ref{eq:odd2}) close to 1/2, and the  expected output of the test will cross these values infinitely often. 
Therefore, we have shown that the expected output of the test $D$ diverges on two independent runs of the process $\rho$, contradicting the consistency of $D$.
This contradiction concludes the proof.	


\begin{thebibliography}{}

\bibitem{nob1} Adams, T.M.,  Nobel, A.B. (1998). 
On density estimation from ergodic processes, {\em The Annals of Probability}
 vol. 26, no. 2, pp. 794--804.

\bibitem{cs} I.~Csiszar, P.~Shields. Information Theory And Statistics: A Tutorial. Now Publishers, 2004.

\bibitem{gml} Gyorfi L., Morvai G., Yakowitz S. (1998), Limits to consistent on-line forecasting for ergodic time series, 
{\em IEEE Transactions on Information Theory} vol. 44 , no. 2, pp. 886--892.


\bibitem{gray}  Gray R. Probability, Random Processes, and Ergodic
Properties. Springer Verlag, 1988.

\bibitem{leh}  Lehmann, E.~L. (1986). Testing Statistical Hypotheses. Springer.


\bibitem{mor2} Morvai G., Weiss B. (2005)  On classifying
 processes. {\em Bernoulli}, vol. 11, no. 3, pp. 523--532.

\bibitem{or2}  D.~S.~Ornstein (1973).  An Application of Ergodic Theory to Probability Theory. 
 {\em The Annals of Probability,} vol.~1, no.~1 pp.~43--58. 

\bibitem{orns} Ornstein, D.~S. (1974).
  Ergodic theory, randomness, and dynamical systems.  Yale
Mathematical Monographs 5, Yale Univ. Press, New Haven, CT.


\bibitem{os} Ornstein, D.~S.,  Shields, P.(1994). The $\bar d$-recognition of processes,  {\em Advances in Mathematics,} vol. 104, pp. 182--224.


\bibitem{ow} 
Ornstein, D.~S. and Weiss, B.(1990). How Sampling Reveals a Process.
 {\em Annals of Probability,} vol. 18 no. 3, pp. 905--930. 

\bibitem{br} Ryabko, B.(1988). Prediction of random sequences and universal coding.
{\em Problems of Information Transmission,} vol. 24, pp. 87--96.


\bibitem{rg} Ryabko, B.,  Astola, J.,  Gammerman, A. (2006).
  Application of Kolmogorov complexity and universal codes to identity testing and nonparametric testing of serial independence for time series.
  {\em Theoretical Computer Science,} vol. 359, pp. 440--448.


\bibitem{rr} Ryabko, D.,  Ryabko, B. (2008). On Hypotheses Testing for Ergodic Processes. 
In Proceedings of IEEE Information Theory Workshop (ITW'08), Porto, Portugal, pp. 281--283. 


\bibitem{sh2} Shields, P.(1993). Two divergence-rate counterexamples, {\em Journal of Theoretical Probability,} vol. 6, pp. 521--545. 

\bibitem{sh}Shields, P.(1998).  The Interactions Between Ergodic
Theory and Information Theory. {\em IEEE Transactions on Information
Theory,} vol. 44, no. 6, pp.~2079--2093.

\bibitem{shbook} Shields, P. (1996)  The Ergodic Theory of Discrete Sample Paths.
AMS Bookstore.


\bibitem{shir} Shiryaev, A. (1996).
\newblock  Probability, second edition.
\newblock Springer.


\end{thebibliography}
\end{document}